\documentclass[numbook, envcountsect, envcountsame, envcountreset, runningheads, smallextended,twocolumn]{svjour3}
\usepackage{graphicx}
\usepackage{fix-cm}
\usepackage{amssymb}
\usepackage{amsmath}
\usepackage{comment, epsfig}
\usepackage{bm}
\usepackage{multirow}
\usepackage{natbib}
\setcitestyle{aysep={}}

\DeclareMathAlphabet{\mathpzc}{OT1}{pzc}{m}{it}

\newtheorem{thm}{Theorem}[section]
\numberwithin{equation}{section}
\allowdisplaybreaks

\newtheorem{cor}{Corollary}[section]

\usepackage[colorlinks, citecolor=blue, linkcolor=red, urlcolor=red]{hyperref}

\begin{document}

\title{Parameter estimation for fractional birth and  fractional death processes}


\author{Dexter O. Cahoy       \and      Federico Polito}


\institute{Dexter O. Cahoy $( \boxtimes )$ \at
       Program of Mathematics and Statistics\\
       College of Engineering and Science\\
		Louisiana Tech University, USA\\
        Tel: +1 318 257 3529\\
        Fax: +1 318 257 2182\\
         \email{dcahoy@latech.edu}
\and
Federico Polito \at
Dipartimento di Matematica\\
 Universit\`a degli Studi di Torino, Italy\\
 Tel: +39 0116702937\\
 Fax: +39 0116702878\\
 \email{federico.polito@unito.it}
}

\date{Received: 08/05/2011 / Re-revised: 09/19/2012}

\maketitle

\begin{abstract}
         The fractional birth and the fractional death processes are more desirable in practice than their classical counterparts as they naturally provide greater flexibility in modeling growing and decreasing systems. In this paper, we propose formal parameter estimation procedures for the fractional Yule, the fractional linear death, and the fractional sublinear  death processes. The methods use all available data possible, are computationally simple and asymptotically unbiased. The  procedures exploited the natural structure of the random inter-birth and inter-death times that are known to be independent but are not identically distributed.  We also showed how these methods can be applied to certain models with more general birth and death rates.  The computational tests showed favorable results for our proposed methods even with relatively small sample sizes. The proposed methods are also illustrated using the branching times of the plethodontid salamanders data of \cite{hal79}.

\end{abstract}

\vspace{-0.2in}
\section{Introduction}

Recently, generalizations of the classical birth and death processes  have been developed using the techniques of fractional calculus. These are called the fractional birth \citep{ucs08,oap10,cap12}  and the fractional death \citep{ops10}  processes, correspondingly. A major advantage of these models over their classical counterparts is that they can capture both Markovian  and non-Markovian structures of a growing or decreasing system.

When the birth and death rates are both linear, they are then called the fractional linear birth or  fractional Yule or Yule--Furry process (fYp) and  fractional linear death process, respectively.  The classical linear birth or Yule process has been  widely used to model various stochastic  systems such as cosmic showers in physics and epidemics in biology to name a few \citep[see e.g., ][]{net94a,dja01, nee01, pe12}. Note  also that the fractional linear birth process was partially investigated by  \cite{ucs08} using the Riemann-Liouville derivative operator but was continued and generalized by \cite{oap10} using the Caputo derivative. The inter-birth time distribution, which provided a way to simulate the  fYp was derived in \cite{cap12}.  With this, we adopt the fYp from  \cite{oap10}. In addition,  the definition of the fractional linear and fractional  sublinear death processes are taken from  \cite{ops10}.   

For completeness, we first enumerate some properties of the fractional Yule (with one progenitor) and the fractional  linear death (with initial population size $n_0 >1$) processes,  which will be used in the subsequent discussions.  Table \ref{t0} below  shows the probability $\tilde{P}(t)$ of no event (no birth or no death) at time $t$,  the state probability mass function $P_i(t)$  or the probability of having  $(i-1)$ births or $(n_0-i)$ deaths by time $t$,  the probability density function $f_i(t)$ of the independent but non-identically distributed random inter-event times, the mean, and the variance of the fractional Yule and the fractional linear death   processes.  Note that the fractional Yule  and  the fractional linear death processes  have the parameters $\lambda>0$ and $\mu>0$ as the birth and death intensities, correspondingly.  


\begin{center}
\begin{table*}[h!t!b!p!]
\caption{Known properties of fractional Yule ($N^\nu(t)$) and linear death ($M^\nu(t)$) processes.} \label{yd}
\hfill{}
\centerline{
\begin{tabular}{|c||c|c|}
\hline
& & \\
& fractional Yule  process & fractional linear death process \\
\hline \hline
& & \\
$\tilde{P}(t)$& $E_{\nu,1}(-\lambda t^{\nu} )$ & $E_{\nu,1}(-\mu n_0 t^{\nu} )$\\
& & \\
$P_i(t)$& $\sum_{j=1}^i \binom{i-1}{j-1} (-1)^{j-1}E_{\nu,1} (-\lambda j t^\nu), \quad i \geq1$ & $\binom{n_0}{i} \sum_{j=0}^{n_0-i} \binom{n_0-i}{j} (-1)^j	E_{\nu,1} (-(i+j) \mu t^\nu), \quad 0 \leq i \leq n_0$\\
& &\\
$f_i(t)$&$ \lambda i t^{\nu-1} E_{\nu,\nu} (-\lambda i t^\nu),\quad i \geq 1 $&  $ \mu(n_0-i) t^{\nu-1} E_{\nu,\nu}(-\mu(n_0-i)t^\nu), \quad 0 \leq i \leq n_0-1 \: (**)$ \\
& & \\
Mean & $E_{\nu,1} \left( \lambda t^\nu \right)$ &$ n_0 E_{\nu,1} \left( - \mu t^\nu \right)$\\
& &\\
Variance & $ 2 E_{\nu,1} \left( 2 \lambda t^\nu \right) - E_{\nu,1} \left( \lambda t^\nu \right) - \left( E_{\nu,1} \left( \lambda t^\nu \right) \right)^2$ &  $ n_0(n_0-1) E_{\nu,1}(-2\mu t^\nu) + n_0 E_{\nu,1}(-\mu t^\nu) - n_0^2 \left(E_{\nu,1}(-\mu t^\nu)\right)^2 \: (**)$ \\
\hline
\end{tabular}
}
\label{t0}
\hfill{}

\end{table*}
\end{center}
Note that
\begin{equation}
			\label{mittagintro}
			E_{\delta,\beta} \left( x \right) = \sum_{j=0}^\infty
			\frac{x^j}{\Gamma \left( \delta j+\beta \right)}
\end{equation}
is the Mittag--Leffler function. \footnote{Note: The entries with  (**) are new results and are derived in Section 2.}

In this article, we propose regression-based procedures to estimate the parameters of the fractional linear birth, the  fractional linear death, and the fractional sublinear death processes. The rest of the paper is organized as follows. In Section 2,  the  specific functional  forms of the inter-death time distributions and the variances of the fractional linear and sublinear death processes are obtained. These results allowed us to apply our methods to these processes. Section 3 introduces the proposed method using the fractional Yule, the fractional linear death, and the fractional sublinear death processes as  examples.  The section also shows some extensions of the procedures to certain models.  Section 4 contains the  empirical test results and the real-data application  of the proposed methods for the case of the fYp only as similar inference procedures can be applied to the fractional linear and fractional sublinear death processes. The summary and  extensions of our study are given in Section 5.

\section{More properties of the fractional linear and  fractional sublinear death processes }

We now derive some properties which will permit us to apply  the proposed estimation procedures to the fractional linear death and the fractional sublinear death processes. More specifically, the theorems below showed that the inter-death times  for  both the fractional linear and sublinear death processes are Mittag-Leffler distributed. The variances of both processes are also derived.

\begin{thm}
	The inter-death time $T_k^\nu$  of the fractional linear  death process $\{M^\nu(t), t>0\}$ with death rate intensity $\mu>0$, and $n_0 \in \mathbb{N}$ initial
	individuals are independent but are non-identically distributed with probability density function\
	\[ 
		\Pr \{T_k^\nu \in dt\}/dt = \mu (n_0-k) t^{\nu-1} E_{\nu,\nu}(-\mu (n_0-k) t^\nu), \notag 
	\]
	where $k=0, 1,  \ldots, n_0-1,$ and  $T_k^\nu$ is the random time separating the $k$th and $(k+1)$th death.
	
	\begin{proof}	
		We prove the theorem by induction. When $k=0$ we obtain
		\begin{align}
			\Pr \{ T_0^\nu \leq t \} & = \Pr \{ M^\nu (t) < n_0\} \\
			& =	1-\Pr \{ M^\nu (t) = n_0\}\notag \\
			& =	1-\tilde{P} (t) \quad (\text{see Table \ref{yd}}) \notag \\
			& = 1 - E_{\nu,1} (-\mu n_0 t^\nu) \notag.
		\end{align}
		Therefore
		\begin{align}
			\Pr \{ T_0^\nu \in dt \} / dt & = \frac{d}{dt} \Pr \{ T_0^\nu \leq t \} \\
			& = \mu n_0 t^{\nu-1} E_{\nu,\nu} (-\mu n_0 t^\nu). \notag
		\end{align}
		For $k=1$ we observe
		\begin{align}
			\Pr &\{ T_0^\nu+T_1^\nu \in dt \} /dt \\
			= {} & \frac{d}{dt} \Pr \{ T_0^\nu+T_1^\nu < t  \} \notag \\
                                   = {} &  \frac{d}{dt}   \Pr \{ M^\nu (t) < n_0 -1 \}  \notag \\
                                   = {} &   \frac{d}{dt} \big[ 1-\Pr \{ M^\nu (t) = n_0  \}  - \Pr \{ M^\nu (t) =  n_0  -1 \} \big]  \notag.
		\end{align}
Using Table \ref{yd} we get
 		\begin{align}
			\Pr &\{ T_0^\nu+T_1^\nu \in dt \} / dt\\
			= {} & -\frac{d}{dt} E_{\nu,1}(-\mu n_0 t^\nu) \notag \\
			& - \frac{d}{dt} \left[ n_0 E_{\nu,1}
			(-(n_0-1)\mu t^\nu) -n_0 E_{\nu,1}(-n_0 \mu t^\nu) \right] \notag \\
			= {} & \mu n_0 t^{\nu-1} E_{\nu,\nu}(-\mu n_0 t^\nu) \notag \\
			& + n_0(n_0-1) \mu t^{\nu-1} E_{\nu,\nu}
			(-\mu(n_0-1)t^\nu) \notag \\
			& - n_0^2 \mu t^{\nu-1} E_{\nu,\nu} (-\mu n_0 t^\nu) \notag \\
			= {} & n_0(n_0-1)\mu t^{\nu-1} \left[ E_{\nu,\nu}(-\mu(n_0-1)t^\nu) \right. \notag \\
			& \left. - E_{\nu,\nu} (-\mu n_0 t^\nu) \right] \notag.
		\end{align}
		To check the preceding results, we can obtain the Laplace transform as
		\begin{align}
			& \int_0^\infty e^{-wt} \Pr \{ T_0^\nu + T_1^\nu \in dt \} \\
			& = \frac{n_0(n_0-1)\mu}{w^\nu + \mu (n_0-1)} - \frac{n_0(n_0-1)\mu}{w^\nu+\mu n_0} \notag \\
			& = \frac{\mu n_0}{w^\nu + \mu n_0} \cdot \frac{\mu(n_0-1)}{w^\nu+\mu(n_0-1)} \notag \\
                        & = \int_0^\infty e^{-ws} \Pr \{ T_0^\nu \in ds \} \int_0^\infty e^{-w y} \Pr \{ T_1^\nu \in ds \}
			\notag \\
			& = \int_0^\infty e^{-wt} \int_0^t \Pr \{ T_1^\nu \in d(t-s) \} \Pr \{ T_0^\nu \in ds \} \notag \\
			& = \int_0^\infty \Pr \{ T_0^\nu \in ds \} \int_s^\infty e^{-zt} \Pr \{ T_1^\nu \in d(t-s) \}, \notag 
		\end{align}
which is just a convolution of two independent variables $T_0^\nu $ and $ T_1^\nu $. For a general $k$ it is sufficient to note that
		\begin{align}
			& \Pr \{ T_0^\nu + \dots + T_k^\nu \in dt \} \\
			& = \int_0^t \Pr \{ T_k^\nu \in d(t-s) \} \Pr \{ T_0^\nu + \dots + T_{k-1}^\nu \in ds \} \notag .
		\end{align}
		By exploiting again the Laplace transform and writing $\mathfrak{D}_k^\nu = T_0^\nu + \dots + T_k^\nu$, we have
		\begin{align}
			& \int_0^\infty e^{-w t} \Pr \{ \mathfrak{D}_k^\nu \in dt \} \\
			& = \int_0^\infty e^{-w t} \int_0^t \Pr \{ T_k^\nu \in d(t-s) \} \Pr \{ \mathfrak{D}_{k-1}^\nu \in ds \}
			\notag \\
			& = \int_0^\infty \Pr \{ \mathfrak{D}_{k-1}^\nu \in ds \} \int_s^\infty e^{-zt} \Pr \{ T_k^\nu \in d(t-s) \}
			\notag \\
			& = \int_0^\infty e^{-w s} \Pr \{ \mathfrak{D}_{k-1}^\nu \in ds \} \int_0^\infty e^{-w y} \Pr \{ T_k^\nu
			\in dy \}\notag \\
			& = \prod_{j=0}^k \int_0^\infty e^{-w s} \Pr \{ T_j^\nu \in ds \} \notag \\
			& = \prod_{j=0}^k \frac{\mu(n_0-j)}{w^\nu + \mu(n_0-j)}. \quad \blacksquare \notag
		\end{align}
	\end{proof}
\end{thm}

We now determine the variance of the fractional linear death process $\{M^\nu(t), t>0\}$. Consider equation (1.6) of \cite{ops10}. That is,
	\begin{align}
		\begin{cases}
			\frac{d^\nu}{dt^\nu} p_k^\nu(t) = \mu (k+1) p_{k+1}^\nu (t) - \mu k p_k^\nu(t), & 0 \leq k \leq n_0, \\
			p_k^\nu(0) =
			\begin{cases}
				1, & k=n_0, \\
				0, & 0 \leq k < n_0.
			\end{cases}
		\end{cases}
	\end{align}
	It is then straightforward to arrive at
	\begin{align}
		\begin{cases}
			\frac{\partial^\nu}{\partial t^\nu} G^\nu(u,t) = -\mu (u-1) \frac{\partial}{\partial u} G^\nu(u,t), \\
			G^\nu(u,0) = u^{n_0},
		\end{cases}
	\end{align}
	where $G^\nu(u,t) = \sum_{k=0}^{n_0} u^k p_k^\nu(t)$ is the probability generating function of the
	fractional linear death process. This in turn leads to
	\begin{align}
		\label{casa}
		\begin{cases}
			\frac{\partial^\nu}{\partial t^\nu} H(t) = -2\mu H(t), \\
			H(0) = n_0(n_0-1),
		\end{cases}
	\end{align}
	where $H(t) = \mathbb{E}(M^\nu(t) (M^\nu(t)-1))$ is the second factorial moment. The solution to \eqref{casa}
	reads
	\begin{align}
		H(t) = n_0(n_0-1) E_{\nu,1}(-2\mu t^\nu),
	\end{align}
	and the variance can be immediately obtained as
	\begin{align}
		\mathbb{V}ar M^\nu(t) & = H(t) + \mathbb{E}M^\nu(t) - (\mathbb{E}M^\nu(t))^2 \\
		& = n_0(n_0-1) E_{\nu,1}(-2\mu t^\nu) \notag \\
        &  + n_0 E_{\nu,1}(-\mu t^\nu) - n_0^2 (E_{\nu,1}(-\mu t^\nu))^2. \notag
	\end{align}
	Note that the above expression, when $\nu=1$, simplifies to the variance of the classical linear
	death process, i.e.\
	\begin{align}
		\mathbb{V}ar M^1(t) = n_0 e^{-\mu t} (1-e^{-\mu t}).
	\end{align}

\begin{figure}[h!t!b!p!]
			\centering
			\includegraphics[height=4.5in, width=3.4in]{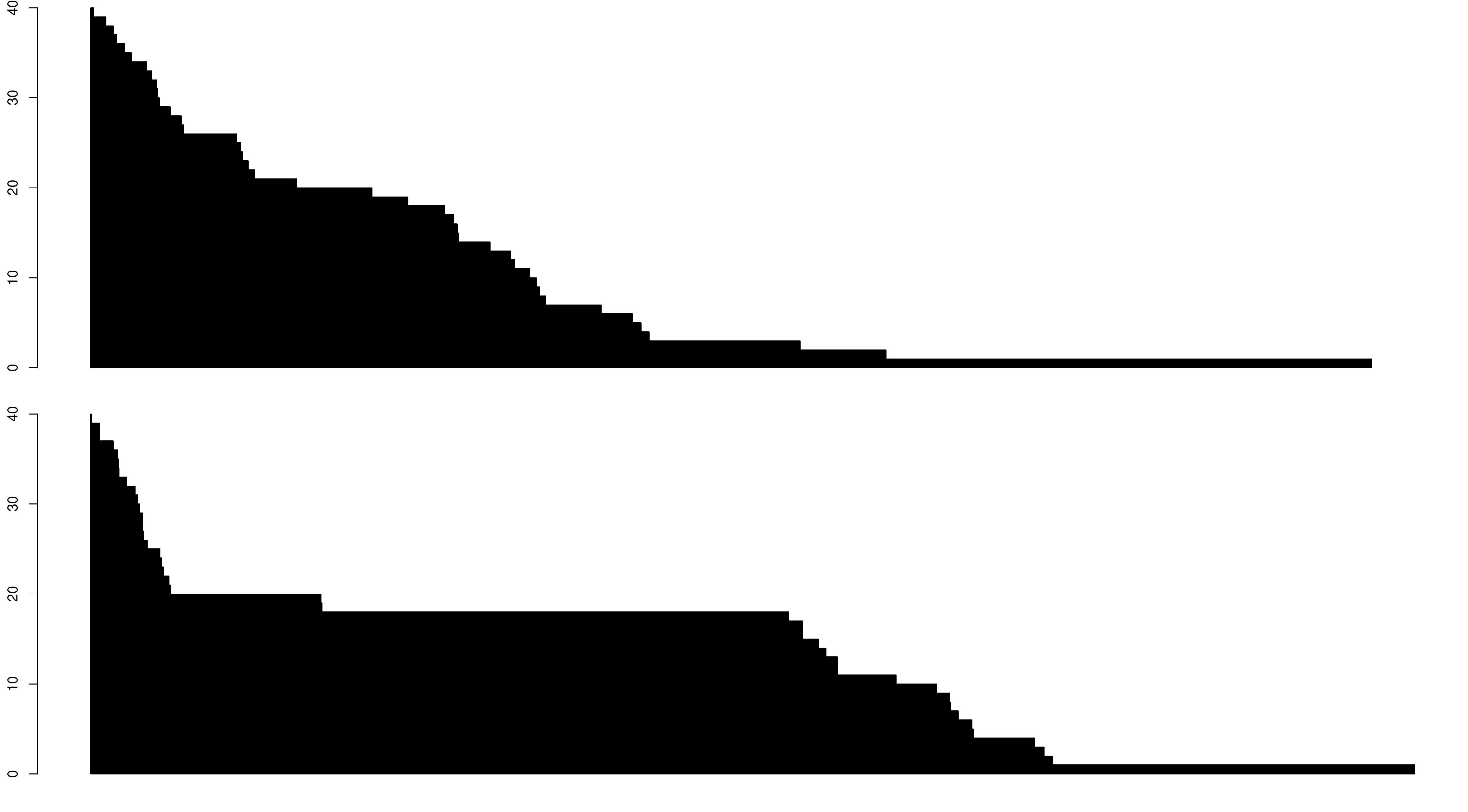}
			\caption{\emph{Sample paths of the classical linear death process  (top) and  the fractional linear death
				process (bottom) in the interval with parameters   $ (\nu, \lambda)=(0.75, 1)$ and initial population size $n_0=40.$  }}
			 \renewcommand\belowcaptionskip{0pt}
			\label{Figure1}
		\end{figure}

Below is the algorithm to generate a typical sample path of a fractional linear  death process in Figure \ref{Figure1}. Note that there are several sub-algorithms  to generate the inter-death times $T_j^\nu$'s that are available in the literature \citep[see e.g.,][]{cap12}.

\vspace{0.1in}
\noindent \textbf{Algorithm:}
		\begin{description}
			\item Step 1.  Let $k=0$ and the population size equal $n_0$.
			\item Step 2. Simulate $T_k^\nu$, and let the $k$th death time be $\mathfrak{D}_k^\nu =
				T_0^\nu + T_1^\nu + T_2^\nu + \cdots +  T_k^\nu.$
			\item Step 3. Set the population size $n_0-k$, and $k=k+1$.
			\item Step 4. Repeat Steps 2--3  for $k=1,\ldots,n_0-1$.
		\end{description}

 It can be gleaned from Figure \ref{Figure1} that the sample path of the fractional linear death process (bottom) seems to decay faster at small times but is slower for large times  than its classical counterpart. The figure also indicates that it is capable of  producing death bursts especially at early stages (corresponding to  small times).

The inter-death time distribution for the fractional sublinear death process  can be easily deduced  (whose proof follows from the previous result and is omitted) from the preceding theorem  as follows.
\begin{thm}
	The fractional sublinear death process $\{\mathfrak{M}^\nu(t)$, $t>0\},$  with death intensity rate $\mu>0$, and $n_0 \in \mathbb{N}$ initial 	 individuals has the following probability density function of the  inter-death times $\mathfrak{T}_k^\nu$'s
	\[
       \Pr \{\mathfrak{T}_k^\nu \in dt\} /dt= \mu (k + 1)  t^{\nu-1} E_{\nu,\nu}(-\mu (k +1) t^\nu), \notag 
	\]
	with $k=0, 1, \ldots, n_0-1$,  where $\mathfrak{T}_k^\nu$ is the random time separating the $k$th and $(k+1)$th death.
\end{thm}

The variance of the fractional sublinear death process  can be determined by considering equation (3.45) of \cite{ops10}. Recall that
	\begin{align}
		\left. \frac{\partial^2}{\partial u^2} \mathfrak{G}^\nu (u,t) \right|_{u=1} & =  \mathbb{E} \left[
		\mathfrak{M}^\nu(t) \left( \mathfrak{M}^\nu(t) - 1 \right) \right] \\
            & = H(t). \notag
	\end{align}
Then
\[
\frac{d^\nu}{d t^\nu} H(t) =   -2 \mu (n_0+1) \left( \mathbb{E}\mathfrak{M}^\nu(t) +  \Pr \{ \mathfrak{M}^\nu(t) = 0 \}
		-1 \right)
\]
\begin{equation}
\qquad \qquad  + \:2 \mu H(t)
\end{equation}
\[
		= -2 \mu (n_0+1) \left( \sum_{k=1}^{n_0} \binom{n_0}{k} (-1)^k E_{\nu,1}(-k\mu t^\nu) \right.
\]
\[
	 \left. +	\sum_{k=1}^{n_0} \binom{n_0+1}{k+1}(-1)^{k+1} E_{\nu,1}(-\mu k t^\nu) \right) + 2\mu H(t)
\]
\[
		= - 2\mu (n_0+1) \left[ \sum_{k=1}^{n_0} \left[ \binom{n_0}{k} - \binom{n_0+1}{k+1} \right]
		(-1)^k E_{\nu,1} (-k\mu t^\nu) \right]
\]
\[
\qquad  \qquad +\: 2\mu H(t)
\]

\[
		=  2 \mu (n_0+1) \sum_{k=1}^{n_0} \binom{n_0}{k+1} (-1)^k E_{\nu,1}(-k\mu t^\nu)
\]
\[
\qquad \qquad +\: 2 \mu H(t).
\]
Using the initial condition $H(0) = n_0(n_0-1)$ and letting $\tilde{H}(w)$ be the Laplace transform of $H(t)$, we write
	\begin{equation}
		w^\nu \tilde{H}(w) - w^{\nu-1} n_0(n_0-1)
	\end{equation} 	
\[
 = 2 \mu (n_0+1) \sum_{k=1}^{n_0} \binom{n_0}{k+1}	(-1)^k \frac{w^{\nu-1}}{w^\nu+k\mu} + 2\mu \tilde{H} (w).
\]
Hence,
\begin{equation}
		\tilde{H}(w)
\end{equation}
\[
 = n_0(n_0-1) \frac{w^{\nu-1}}{w^\nu-2\mu} + 2\mu (n_0+1) \sum_{k=1}^{n_0}
		\binom{n_0}{k+1} \times
\]
\[
	\qquad (-1)^k w^{\nu-1} \frac{1}{(w^\nu + k \mu)(w^\nu-2\mu)}
\]
\[
= n_0(n_0-1) \frac{w^{\nu-1}}{w^\nu-2\mu} + 2\mu (n_0+1)
\]
\[
\times \sum_{k=1}^{n_0}
		\binom{n_0}{k+1} (-1)^k w^{\nu-1} \left[ \frac{1}{w^\nu+k\mu} - \frac{1}{w^\nu-2\mu} \right]
		\frac{1}{(-2\mu)}
\]
\[
= n_0(n_0-1) \frac{w^{\nu-1}}{w^\nu-2\mu} + \frac{w^{\nu-1}}{w^\nu-2\mu} (n_0+1)
		\sum_{k=1}^{n_0} \binom{n_0}{k+1} (-1)^k
\]
\[
	\qquad - (n_0+1) \sum_{k=1}^{n_0} \binom{n_0}{k+1} (-1)^k \frac{w^{\nu-1}}{w^\nu+k\mu}
 \]
 \[
		= \frac{w^{\nu-1}}{w^\nu-2\mu} (1-n_0) - (n_0+1) \sum_{k=1}^{n_0} \binom{n_0}{k+1} (-1)^k
		\frac{w^{\nu-1}}{w^\nu+k \mu}.
\]
The second factorial moment can be easily shown as
\begin{equation}
		H(t) =  - (n_0-1) E_{\nu,1}(2\mu t^\nu)
\end{equation}
\[
+ (n_0+1) \sum_{k=1}^{n_0} \binom{n_0}{k+1}
		(-1)^{k+1} E_{\nu,1}(-k\mu t^\nu).
\]
Thus, the variance simply follows as
\[
		\mathbb{V}ar \, \mathfrak{M}^\nu(t) = H(t) + \mathbb{E} \mathfrak{M}^\nu(t) -
		\left[ \mathfrak{M}^\nu(t) \right]^2
\]
\[		
= - (n_0-1) E_{\nu,1}(2\mu t^\nu) + (n_0+1) \sum_{k=1}^{n_0} \binom{n_0}{k+1}
\]
\[
\qquad \qquad \times	(-1)^{k+1} E_{\nu,1}(-k\mu t^\nu)
\]
\[
\qquad + \sum_{k=1}^{n_0} \binom{n_0+1}{k+1} (-1)^{k+1} E_{\nu,1}(-\mu k t^\nu)
\]
\[	
\qquad - \left[ \sum_{k=1}^{n_0} \binom{n_0+1}{k+1} (-1)^{k+1} E_{\nu,1}(-k \mu t^\nu) \right]^2.
\]
Note that the algorithm above could be easily adopted to simulate sample trajectories of the fractional sublinear death process. 

\section{Parameter estimation}

\subsection{Estimation for the fractional Yule or linear birth process}

We  now illustrate our estimation approach for the fYp with birth rate $\lambda i, i \geq 1.$  Furthermore,  assume that a sample trajectory of $n$ births corresponding to $n$ random inter-birth times $T_i$'s of the fractional linear birth process is observed.  That is,  $n$ independent but are not identically distributed random inter-birth times of the fractional linear birth process are given.  This also insinuates that  only  a single datum is obtained from each of the $n$ different  Mittag-Leffler distributions. This observation  and the Mittag-Leffler's seemingly complex structure pose a computational challenge on how to estimate the model parameters more efficiently especially for small population sizes.  Recall the structural representation  of the Mittag-Leffler distributed random inter-birth time $T_i  \stackrel{d}{=}  E^{1 / \nu} S_\nu$  \citep[see][]{cap12},   where $E \stackrel{d}{=} \exp ( \lambda i)$  is independent of $S_\nu$  which is a one-sided $\alpha^+$-stable distributed random variable. Applying  the logarithmic transformation and taking  the expectation on both sides,  it  can  be   easily shown   that  the mean and variance \citep[see details in][]{cuw10} of the log-transformed $i$-th  random  sojourn time $T_i^{'}=\ln \left( T_i \right)$  of the fYp  are
\begin{equation} \label{assume1}
\mu_{T_i^{'}}= \frac{-\ln \left( \lambda i \right) }{\nu} - \gamma,
\end{equation}
and
\begin{equation}
 \sigma_{T_i^{'}}^2 =\pi^2 \left( \frac{1}{3\nu^2} -\frac{1}{6} \right),
\end{equation}
respectively, where $ \gamma \approx 0.5772156649$ is the  Euler - Mascheroni's constant. The first two moments above therefore suggest that the following simple linear regression model can be fitted/formulated:
\begin{equation}
T_i^{'} = a_0 + a_1 \ln i + \varepsilon_i, \qquad  i =1,\ldots, n, 
\end{equation}
where
\begin{equation}
\label{b1e}
a_0= \frac{-\ln (\lambda ) }{\nu} - \gamma,\qquad a_1= \frac{-1}{\nu},
\end{equation}
and $\varepsilon_i \stackrel{iid}{=} \textsl{N} \left(\mu_\varepsilon=0,  \sigma_\varepsilon^2 = \sigma_{T_i^{'}}^2 \right)$. The trick used here was to factor out the non-identical means of the log-transformed random inter-birth or sojourn times, which are linear functions of the logarithm of the   known fixed $i$. Thus, this leads to studying the widely used simple  linear regression model \citep[see][]{mpv06}.

\subsubsection{Point estimation}

Inverting the least squares (LS) estimators 
\begin{equation}
 \hat{a}_1 = \frac{ \sum_{j=1}^n  T_j^{'} \left( \ln j -   \overline{\ln i } \right) }{\sum_{j=1}^n  \left( \ln j -   \overline{\ln i } \right)^2}
\end{equation} 
and $\hat{a}_0= \overline{T_i^{'} }  -  \hat{a}_1 \cdot \overline{\ln i }$ gives the LS-based point estimators of $\nu$ and $\lambda$ as
\begin{equation}
\label{nulsp}
\widehat{\nu}_{ls}= \frac{-1}{\widehat{a}_1}
\end{equation}
and
\begin{equation}
\label{mulsp}
\widehat{\lambda}_{ls}= \exp\left( \left( \widehat{a}_0 + \gamma  \right) \big/ \; \widehat{a}_1 \right),
\end{equation}
respectively, where  $\overline{\ln i } = \sum\limits_{j=1}^n \ln j /n$, and $ \overline{T^{'} }=\sum\limits_{j=1}^n T_j^{'} /n$.  Equating $\sigma_\varepsilon^2 $ or  $\sigma_{T_i^{'}}^2$ in  (\ref{assume1})   with its unbiased estimator
\begin{align}
\widehat{\sigma}_u^2   = \sum_{j=1}^n   \widehat{\varepsilon}_j^2/ (n-2),
\end{align}
we get the residual-based point estimators
\begin{equation}
\label{nuresp}
\widehat{\nu}_{res}= \frac{1}{\sqrt{3 \left( \widehat{\sigma}_u^2 \big/ \pi^2 +\frac{1}{6} \right)} } \quad \text{\citep[see][]{cuw10}}
\end{equation}
and
\begin{equation}
\label{muresp}
\widehat{\lambda}_{res}= \exp\left( -\widehat{\nu}_{res} \left( \widehat{a}_0 + \gamma  \right) \right)
\end{equation}
of the model parameters $\nu$  and $\lambda$, correspondingly where $ \widehat{\varepsilon}_i  = T_i^{'} - \widehat{T_i^{'}}$, and $ \widehat{T_i^{'} }=\widehat{}a_0 + \widehat{a}_1 \ln i $. Note that  the residual-based estimators exploit the residuals to estimate $\nu$ rather than the negative inverse of the LS estimate of the slope $a_1$.

\subsubsection{Interval estimation}


We now develop interval estimators using the large-sample properties of the least squares estimators   $\widehat{b}_0$ and $\widehat{b}_1$ above. The following result shows the joint asymptotic behavior of the proposed point estimators of $\nu$ and $\lambda$ for the fYp.
\begin{thm}
			\label{thea}
			Let $0<\nu \leq 1$ and $\lambda >0$.  Then
\[
\sqrt{n}\left(
  \begin{array}{c}
    \widehat{\nu}_{ls} - \nu  \\
    \widehat{\lambda}_{ls}-\lambda \\
  \end{array}
\right) \stackrel{d}{\longrightarrow}  \textsl{N} \left[\bm{0} , n  \sigma_\varepsilon^2 \bf{C}  \right]
\]
where ``$ \stackrel{d}{\longrightarrow}$" denotes convergence in distribution,
\[
\bf{C} = \left(
       \begin{array}{ccccc}
        C_1  &  & C_{12} \\
         &  &   \\
      C_{21} &  &   C_2  \\
       \end{array}
     \right),
\]
$C_1=\nu^4  s^{-1}$,\\
$C_{12}=C_{21}=\lambda  \nu^3 \left(  \left( \overline{\ln i }  + \ln (\lambda) \right)/ s \right)$,\\
$C_2= \left(\nu \lambda \right)^2 \left( 1/n +  \left( \overline{\ln i }^2\right. \right.$\\
$ \left. \left. \qquad + 2\ln (\lambda) \overline{\ln i } + ( \ln (\lambda) )^2 \right) / s  \right)$,\\ and $s=\sum\limits_{j=1}^n \left(\ln j - \overline{\ln i }\right)^2$. 

\begin{proof}
Recall the large-sample normality of the least\\ squares estimators  $\widehat{a}_0$ and $\widehat{a}_1$, i.e.,
\[
\sqrt{n}\left(
  \begin{array}{c}
    \widehat{a}_0 - a_0 \\
    \widehat{a}_1 - a_1  \\
  \end{array}
\right) \stackrel{d}{\longrightarrow}  \textsl{N} \left[\bm{0} ,  \bf{\Sigma} \right]
\]
 where the covariance matrix $\bf{\Sigma}$ is defined as
\[
\bf{\Sigma} = \text{n} \sigma_\varepsilon^2 \left(
       \begin{array}{cccc}
        \left( 1/n + \overline{\ln i }^2 \big/s \right)  & & & -\overline{\ln i }/s\\
         &  &  & \\
       - \overline{\ln i }/s  & &  &  s^{-1} \\
       \end{array}
     \right).
\]
Recall the multivariate delta method \citep{fer96}:   If $\sqrt{n}\big(\widehat{\bm{\beta}}_n- \bm{\beta}\big) \longrightarrow  \textsl{N} \left[\bm{0} ,  \bf{\Sigma} \right]$  then 
\[
\sqrt{n}\big(\textbf{g}(\widehat{\bm{\beta}}_n)-\textbf{g}(\bm{\beta})\big)\stackrel{d}{\to} \textsl{N}\left[ \bm{0} ,\; \bm{\dot{\textbf{g}}}(\bm{\beta})^{\text{T}}\bf{\Sigma}\bf{\dot{g}}(\bm{\beta})\right].
\]
Hence, using the delta method above, $\widehat{\bm{\beta}}_n=(\widehat{a}_0 , \widehat{a}_1 )^\text{T},  \\ \textbf{g}(\widehat{\bm{\beta}}_n)= (\widehat{\nu}_{ls},\widehat{\lambda}_{ls})^\text{T}$,  $\textbf{g}(\bm{\beta})= (\nu,\lambda)^\text{T}$,  and  the Jacobian matrix
\[
\bm{\dot{\textbf{g}}}(\bm{\beta}) = \left(
       \begin{array}{ccccc}
       0& & &  & \exp\left( (a_0 + \gamma) / a_1 \right)/a_1 \\
         &  &  &  &\\
      1/a_1^2  & &  & &  -\exp\left( (a_0 + \gamma) / a_1 \right) \left( a_0 + \gamma \right)/ a_1^2  \\
       \end{array}
     \right),
\]
we obtain the final expression of the covariance matrix by simply substituting back $a_0 = -\ln (\lambda ) / \nu - \gamma$ and $a_1= -1 / \nu. \quad \blacksquare $
\end{proof}
\end{thm}

\begin{cor}
Approximate  $(1-\alpha)100\%$ confidence intervals for  $\nu$ and $\lambda$ can be deduced as
\begin{equation}
\label{nulsi}
\widehat{\nu}_{ls} \; \; \pm \; \; z_{\alpha/2} \widehat{\sigma}_\varepsilon  \widehat{\nu}_{ls}^2 \sqrt{s^{-1}},
\end{equation}
and
\begin{align}
 \widehat{\lambda}_{ls} &  \; \; \pm \; \; z_{\alpha/2} \widehat{\sigma}_\varepsilon \widehat{\nu}_{ls} \widehat{\lambda}_{ls} \Bigl( 1/n \notag \\ \label{mulsi}
& + \; \; \left( \overline{\ln i }^2 + 2\ln ( \widehat{\lambda}_{ls}) \overline{\ln i } + ( \ln ( \widehat{\lambda}_{ls}) )^2 \right) \big/ s   \Bigr)^{1/2},  
\end{align}
respectively, where $z_{\alpha/2}$ is the $(1-\alpha/2)$th quantile of the standard normal distribution and $0 <\alpha<1$.
\end{cor}		


We now propose another interval estimators  which utilize the residual-based estimate of $\nu$, and a bootstrap technique.    It can be inferred from \cite{cuw10} that a residual-based $(1-\alpha)100\%$ confidence interval for $\nu$ can be
\begin{equation}
\label{nuresi}
\widehat{\nu}_{res} \; \; \pm \; \; z_{\alpha/2}  \sqrt{\frac{\widehat{\nu}_{res}^2\left(32-20\widehat{\nu}_{res} ^2-\widehat{\nu}_{res} ^4\right)}{40 n} },
\end{equation}
where $z_{\alpha/2}$ is defined above.
A residual-based $(1-\alpha)100\%$ interval estimate for $\lambda$   can also be 
\[
\label{muresi}
\widehat{\lambda}_{res} \; \; \pm \; \; z_{\alpha/2}  \Bigg[ \frac{  e^{ -2\widehat{\nu}_{res} (\widehat{a}_0 + \gamma  ) }  \bigg( \frac{\widehat{\nu}_{res}^2\left(32-20\widehat{\nu}_{res} ^2 -\widehat{\nu}_{res} ^4\right) }{40} } {n} 
\]
\begin{equation}
\hspace{0.5in}  + \; \; \widehat{\nu}_{res}^2 \widehat{\sigma}_u^2 \left(  1/n + \overline{ \ln i}^2/s \right) \bigg)  \Bigg]^{1/2}.
\end{equation}
Since the small-sample performance  of  $\widehat{\nu}_{res}$  and the residual-based  interval estimator  in (\ref{nuresi})  have been shown to perform well already \citep[see, e.g.,][]{cuw10},  we apply a non-parametric percentile bootstrap  technique to $\widehat{\lambda}_{res}$ using the fixed-regressor approach  to obtain a small-sample interval estimator of $\lambda$. This well-known procedure  is slightly modified by first dividing each residual $\widehat{\varepsilon}$ by $\sqrt{1-h_i},$  where $h_i$  is the $i$th leverage or the $i$th diagonal entry in the hat matrix  before sampling from the transformed residuals. Note that the division of $\sqrt{1-h_i}$  is simply for correction as  the true variance of the residual $\widehat{\varepsilon}_i $ is  $\mathbb{V}ar \;  \widehat{\varepsilon}_i =  \sigma_\varepsilon^2(1-h_i)$  \citep[see][]{mpv06}.   Hence, the bootstrap counterpart of $\widehat{\lambda}_{res}$ is calculated as $\widehat{\lambda}_{res}^{*}=  \exp\left( -\widehat{\nu}_{res}^{*} \left( \widehat{a}_0^{*} + \gamma  \right) \right)$ where $\widehat{\nu}_{res}^{*}$  used the bootstrapped transformed or weighted residuals. A clear advantage of the asymptotic-based procedures over the re-sampling-based ones is  that they are faster  to calculate especially for large sample sizes.

\subsection{Estimation for the fractional linear and the fractional sublinear death processes}

Assuming that a sample trajectory of $n_0$ deaths  corresponding to $n_0$ random inter-death times $T_k^\nu$'s of a fractional death process is observed. Following the procedure for the fractional linear birth process in the preceding subsection, we can estimate the parameters $\nu$ and  $\mu$  by regressing $\ln \left(T_k^\nu \right)$ with $\ln (n_0-k)$. That is,  we fit the following simple linear regression model:
\begin{equation}
\ln \left(T_k^\nu \right)= b_0 + b_1 \ln (n_0-k) + \varepsilon_k,
\end{equation}
where  $k =0,\ldots, n_0-1,$ 
$b_0=  -\ln (\mu) / \nu - \gamma$, 
$b_1$  is given in  (\ref{b1e}) of subsection 3.1, and $\varepsilon_k \stackrel{iid}{=} \textsl{N} \left(\mu_\varepsilon=0,  \sigma_{\ln \left(T_k^\nu \right)}^2 \right)$. Following the methodology in the preceding subsection,  we can straightforwardly obtain the corresponding LS-based point estimates of $\nu$ and $\mu$ from (\ref{nulsp}) and (\ref{mulsp})  as 
\begin{equation}
\label{nulsp2}
\widehat{\nu}_{ls}= \frac{-1}{\widehat{b}_1}
\end{equation}
and
\begin{equation}
\label{mulsp2}
\widehat{\mu}_{ls}= \exp\left( \left( \widehat{b}_0 + \gamma  \right) \big/ \; \widehat{b}_1 \right),
\end{equation}
respectively,  where 
 $\overline{\ln (n_0- k) } = \sum\limits_{j=0}^{n_0-1}\frac{\ln (n_0-j)}{n_0}$, $ \overline{\ln \left(T^\nu \right) }=\sum\limits_{j=0}^{n_0-1}  \frac{\ln \left(T_j^\nu \right)}{n_0}$, $\hat{b}_0=  \overline{\ln \left(T^\nu \right) } -  \hat{b}_1 \cdot  \overline{\ln (n_0- k)}$, and 
\begin{equation}
 \hat{b}_1 = \frac{ \sum_{j=0}^{n_0-1} \ln \left(T_j^\nu \right)  \left( \ln (n_0-j) - \overline{\ln (n_0- k) }  \right) }{\sum_{j=0}^{n_0-1}  \left( \ln (n_0-j) -   \overline{\ln (n_0- k) }  \right)^2}. 
\end{equation} 
Furthermore,  the LS-based interval estimates for $\nu$ and $\mu$ directly follow from (\ref{nulsi}) and (\ref{mulsi}) of Corollary 3.1 in subsection 3.1.2,  correspondingly.   Hence, the approximate $(1-\alpha)100\%$  for $\nu$ and $\mu$ can be explicitly written as
\[
\label{nulsi2}
\widehat{\nu}_{ls} \;  \pm \; z_{\alpha/2} \widehat{\sigma}_\varepsilon  \widehat{\nu}_{ls}^2 \sqrt{\left( \sum \limits_{j=0}^{n_0-1} \left( \ln (n_0-j)- \overline{\ln (n_0-k)} \right)^2\right)^{-1}},
\]
and
\begin{align}
\label{mulsi2}
 \widehat{\lambda}_{ls} &  \; \; \pm \; \; z_{\alpha/2} \widehat{\sigma}_\varepsilon \widehat{\nu}_{ls} \widehat{\lambda}_{ls} \Bigl( 1/n \notag \\
& \hspace{-0.2in}  +  \left( \overline{\ln (n_0-k) }^2 + 2\ln ( \widehat{\lambda}_{ls}) \overline{\ln (n_0-k) } + ( \ln ( \widehat{\lambda}_{ls}) )^2 \right) \big/ s   \Bigr)^{1/2},
\notag
\end{align}
correspondingly. On the other hand, the residual-based point and interval estimators of $\nu$ and $\mu$ immediately follow from subsection 3.1 as well, where  $\hat{a}_0$  is replaced by $\hat{b}_0$ in (\ref{muresp}), $\widehat{\varepsilon}_k  = \ln \left(T_k^\nu \right) - \widehat{\ln \left(T_k^\nu \right)}$, and $ \widehat{\ln \left(T_k^\nu \right)}=\widehat{b}_0 + \widehat{b}_1 \ln (n_0-k)$.

A similar approach can be done  to obtain estimates  for the fractional  sublinear death process. That is,  we regress $\ln \left(\mathfrak{T}_k^\nu \right)$ with $\ln (k+1), k=0, 1, \ldots, n_0-1,$  or fit the model
\begin{equation}
\ln \left(\mathfrak{T}_k^\nu \right)= c_0 + c_1 \ln (k+1) + \varepsilon_k,
\end{equation}
and follow the procedures used for fractional Yule and the fractional linear death processes.  In general, we simply replace $\lambda, \ln i,  \overline{\ln i}$   by $\mu, \ln (n_0-k)$ or $\ln (k+1)$, and $\overline{\ln (n_0-k)}$  or $\overline{\ln (k+1)}$, accordingly  in the methods of subsection 3.1 to obtain the parameter estimators for the fractional linear death and the fractional sublinear death processes.

\subsection{Some Extensions}

Assume that a fractional birth or death process exists with rates $\theta_j, j=1, 2, \ldots, n$,  where the $j$th  inter-event time $X_j$  is Mittag-Leffler distributed with parameter $\theta_j$.  Then  the mean  of $X_j^{'}=\ln \left( X_j \right)$  is
\begin{equation}
\mu_{X_j^{'}}= \frac{-\ln \left( \theta_j \right) }{\nu} - \gamma.
\end{equation}
Based on the above mean formulation, we use the model 
\begin{equation}
X_j^{'} =  d_0 + d_1 \cdot q(j) + \varepsilon_j
\end{equation}
to estimate more  forms of the parameters or rates under the two cases  below.

\vspace{0.1in}
\noindent \emph{Case 1: When $\ln (\theta_j) = m(\theta) +  q(j)$ for  some appropriate  known functions $m(\theta)$ and $q(j)$ of the parameter $\theta$  and $j \in \mathbb{N}$, correspondingly.}

\vspace{0.1in}

In this case, the general form of the regression model  that could be used for estimation is
\begin{equation}
X_j^{'} = -\left(\gamma + \frac{m(\theta)}{\nu} \right)-  \frac{1}{\nu}\cdot q(j) + \varepsilon_j.
\end{equation}
Clearly, $d_0 = -(\gamma +  m(\theta) / \nu),  d_1=  - 1/\nu$, and  $q(j)$  is the regressor variable.  Using $\widehat{\nu}_{res}$  or $-1/\widehat{d}_1$  and inverting the least squares estimate $ \widehat{b}_0$, we can compute $\widehat{m(\theta)}$ and $\widehat{\theta}$ sequentially.  Note that the explicitness of $\widehat{\theta}$ depends on the form of $m$. 

\vspace{0.1in}

 \emph{Example 1}: When $\theta_j $ is linear, i.e.,  $\theta_j = \theta j$  then  $\ln (\theta_j)= \ln(\theta) + \ln(j)$, where $m(\theta) = \ln (\theta)$ and $q(j)=\ln(j)$, respectively. Note that this parametrization corresponds to the fractional Yule, the fractional linear death, and the fractional sublinear death processes.
 
\vspace{0.1in}

 \emph{Example 2}:  If $\theta_j =  e^{\theta + j}$  then $\ln (\theta_j)=  \theta +  j$, where $m(\theta) = \theta $ and $q(j)=j$, correspondingly. This suggests that $d_0 = - (\gamma + \theta/ \nu)$.

\vspace{0.1in}

\noindent \emph{Case 2: When $\ln (\theta_j) = m(\theta) \cdot q(j)$ for  some appropriate  known functions $m(\theta)$ and $q(j)$ of the parameter $\theta$  and $j \in \mathbb{N}$, correspondingly. }
\vspace{0.1in}

The general form of the regression model in this case is
\begin{equation}
X_j^{'} = -\gamma-  \frac{m(\theta)}{\nu}  \cdot q(j) + \varepsilon_j.
\end{equation}
Apparently, $d_0 = -\gamma, d_1= -m(\theta)/\nu $, and $q(j) $ is the predictor variable.  Using $\widehat{\nu}_{res}$ and inverting the least squares estimate $ \widehat{d}_1$, we can calculate $\widehat{m(\theta)}$ and $\widehat{\theta}$ successively.

\vspace{0.1in}

 \emph{Example 1}: If $\theta_j = \theta^j$  then  $\ln (\theta_j)= \ln(\theta)\cdot j$, where $m(\theta) = \ln (\theta)$ and $q(j)= j$, correspondingly.  This indicates that $d_1 = - \ln (\theta)/ \nu$. 

\vspace{0.1in}

 \emph{Example 2}:  When $\theta_j =  e^{\theta \cdot j}$  then $\ln (\theta_j)=  \theta \cdot  j$, where $m(\theta) = \theta $ and $q(j)=j$, respectively. This shows that $d_1 = - \theta/ \nu$.

\section{Method testing and application}

\subsection{Empirical test}

For the sake of reproducibility, we now test our procedures using the fYp as a particular example as similar approach can be carried out for both the fractional linear death and the fractional sublinear death processes. In point estimation testing, we evaluated the finite-sample properties (unbiasedness and homogeneity) by computing the average and the  median absolute deviation (MAD) of the estimates using 1000 simulations for sample sizes $n=100, 500,$ and $1000$.  These values are shown in Table \ref{t1} below.  The relative fluctuation (RF=100\% $\times$MAD/mean) of $\widehat{\nu}_{ls}$ decreases from 19.23\% (corresponds to $\nu=0.1,\; n=100$) to as little as 4.41\% (with $\nu=0.95$ and $n=1000$).  On the other hand, the residual-based $\widehat{\nu}_{res}$'s RF ranges from  4.41\% (with $\nu=0.95$ and $n=1000$)   to 1.89\% (corresponds to $\nu=0.95,\; n=1000$). While $\widehat{\lambda}_{ls}$'s RF improves from 33.94\% (corresponds to $\lambda=0.5,\; n=100$) to 30.48\% (with $\lambda=5$ and $n=1000$), $\widehat{\lambda}_{res}$'s RF decays faster from  55\% ($\lambda=1,\; n=100$) to 25.29\% ($\lambda=5$ and $n=1000$).  In general, the relative fluctuations of the residual-based estimators tend to decay faster than the LS-based estimates. They are also less bias than the LS-based estimators especially for $n \leq 100$. Nonetheless, both the residual- and LS-based point estimators are asymptotically unbiased as expected.

\begin{center}
\begin{table*}[h!t!b!p!]
\caption{\emph{Mean point estimates of and  dispersions from the true parameters  $\nu$ and $\lambda$.}} \centerline {
\begin{tabular*}{4.5in}{@{\extracolsep{\fill}}c|c@{\hspace{0.1in}}||cc|cc|cc}
 \hline
\multirow{2}{*}{$(\nu, \lambda)$} & \multirow{2}{*}{Estimator}  &  \multicolumn{2}{c|}{$n=100$}  &  \multicolumn{2}{c|}{$n=500$} &  \multicolumn{2}{c}{$n=1000$}  \\
 &   & Mean & MAD& Mean& MAD & Mean & MAD \\
  \hline \hline
\multirow{4}{*}{$(0.1, 1)$}  &  $\widehat{\nu}_{ls}$     & 0.104  & 0.020    &0.101  & 0.008 &  0.100 & 0.006  \\
               &  $\widehat{\nu}_{res}$     & 0.103  & 0.008    &0.100  & 0.004  & 0.100 & 0.003  \\ \cline{2-8}
               &  $\widehat{\lambda}_{ls}$     & 3.190  & 0.665    &1.151  & 0.407 &  1.077 & 0.318  \\
               &  $\widehat{\lambda}_{res}$     & 1.318  & 0.725    &1.091  & 0.408  & 1.051 & 0.322  \\
\hline
\multirow{4}{*}{$(0.25, 0.1)$}  &  $\widehat{\nu}_{ls}$     & 0.261  & 0.048    &0.252  & 0.021 &  0.251 & 0.014  \\
               &  $\widehat{\nu}_{res}$     & 0.252  & 0.022    &0.251  & 0.010  & 0.250 & 0.007  \\ \cline{2-8}
               &  $\widehat{\lambda}_{ls}$     & 0.119  & 0.031    &0.106  & 0.025 &  0.103 & 0.021  \\
               &  $\widehat{\lambda}_{res}$     & 0.131  & 0.071    & 0.109  & 0.044  & 0.106 & 0.033  \\
\hline
\multirow{4}{*}{$(0.5, 0.5)$}  &  $\widehat{\nu}_{ls}$     & 0.521  & 0.100    &0.505  & 0.040 &  0.501 & 0.028  \\
               &  $\widehat{\nu}_{res}$     & 0.506  & 0.041    &0.501  & 0.018  & 0.500 & 0.014  \\ \cline{2-8}
               &  $\widehat{\lambda}_{ls}$     & 0.825  & 0.280    &0.565  & 0.190 &  0.532 & 0.142  \\
               &  $\widehat{\lambda}_{res}$     & 0.640  & 0.350    &0.555  & 0.216  & 0.528 & 0.164  \\
  \hline
\multirow{4}{*}{$(0.75, 0.25)$}  &  $\widehat{\nu}_{ls}$     & 0.774  & 0.121    &0.755  & 0.052 &  0.751 & 0.036  \\
               &  $\widehat{\nu}_{res}$     & 0.755  & 0.056    &0.752  & 0.023  & 0.750 & 0.016  \\ \cline{2-8}
               &  $\widehat{\lambda}_{ls}$     & 0.313  & 0.093    &0.266  & 0.069 &  0.259 & 0.056  \\
               &  $\widehat{\lambda}_{res}$     & 0.300  & 0.146    &0.268 & 0.094  & 0.260 & 0.072  \\
 \hline
\multirow{4}{*}{$(0.95, 5)$}  &  $\widehat{\nu}_{ls}$     & 0.969  & 0.131    &0.953  & 0.058 &  0.952 & 0.042  \\
               &  $\widehat{\nu}_{res}$     & 0.955  & 0.055    &0.950  & 0.024  & 0.950 & 0.018  \\ \cline{2-8}
               &  $\widehat{\lambda}_{ls}$  & 11.251   & 3.492  & 5.836  & 2.104 & 5.397 & 1.645\\
               &  $\widehat{\lambda}_{res}$     & 5.978  & 2.544    &5.375  & 1.635  & 5.206 & 1.317  \\
  \hline
\end{tabular*}
}
\label{t1}
\end{table*}
\end{center}


Table \ref{t2}  below shows the  averaged lower and upper $95\%$ confidence bounds using the formulae in Section 3. These bounds used 1000 simulation runs for each of the sample sizes $n=15, 30, 100,$ and 500. Note that the residual-based interval estimator $\widehat{\lambda}_{res}^{*}$ utilized 500 bootstrap samples and $\widehat{\sigma}_{\varepsilon}^2 = \pi^2 \left( 1/(3 \widehat{\nu}_{res}^2) -1/6 \right)$  is used to estimate the error variance in our LS-based procedures. Observe that  some of the interval estimates for sample sizes $n=15$, and $n=30$ are omitted as they are unreliable  due to the  multiple error warnings that showed up during the computation process.  Moreover, the convergence of the coverage probabilities to their  true levels for the LS-based method  is  made faster by using the error variance estimate $ \widehat{\sigma}_\varepsilon^2= \pi^2 \left( 1/(3 \widehat{\nu}_{ls}^2) -1/6 \right)$.  From Table \ref{t2}, it is apparent that the residual-based interval estimates of $\nu$ are narrower and are  better centered  around the true parameter values than the least-squares' even when the sample size is as large as 500. In addition, our simulations showed that the asymptotic or non-bootstrapped  residual-based interval estimator of $\lambda$ gives more sensible results than the LS-based procedure for small samples.  Nevertheless, the LS-based interval estimates for $\lambda$ are  more accurately centered than the bootstrapped residual-based  estimates especially  for large samples.

\begin{center}
\begin{table*}[h!t!b!p!]
\caption{\emph{Average 95\% confidence intervals  for different values of $\nu$ and  $\lambda$.}}
 \centerline {
\begin{tabular*}{5.5in}{c@{\hspace{0.05in}}|c@{\hspace{0.05in}}||c@{\hspace{0.05in}}|c@{\hspace{0.05in}}|c@{\hspace{0.05in}}|c@{\hspace{0.05in}}}
\hline
$(\nu, \lambda)$ & Estimator &  $n=15 $  & $n=30$ &  $n=100$ &  $n=500$  \\
  \hline \hline
\multirow{4}{*}{$(0.1, 1)$}  &  $\widehat{\nu}_{ls}$   & & & (0.063 ,  0.142) & (0.084 ,  0.117)   \\
                &  $\widehat{\nu}_{res}$    & (0.060 , 0.158)   &(0.071 , 0.137)  & (0.083 , 0.118) & (0.092 ,  0.108) \\ \cline{2-6}
                &  $\widehat{\lambda}_{ls}$    &    &   & (-35.8067 , 58.864) & (0.108 ,  2.228) \\
               &  $\widehat{\lambda}_{res}$    &    &   & (-0.362 ,   2.885) & (0.183 ,  2.038) \\
   &  $\widehat{\lambda}_{res}^{*} $    & (0.214 , 53.963)   &(0.243 , 16.667)  & (0.297 ,  5.560) & (0.464 ,  2.652) \\ 
\hline
\multirow{4}{*}{$(0.25, 0.1)$}  &  $\widehat{\nu}_{ls}$   &   &   & (0.162 ,  0.359)  & (0.211 ,  0.292)   \\
                &  $\widehat{\nu}_{res}$    & (0.151 , 0.389)   &(0.211 , 0.298)  & (0.209 ,  0.296) & (0.231 ,  0.269) \\ \cline{2-6}
                &  $\widehat{\lambda}_{ls}$    &  &  & (0.034 ,  0.201)  & (0.052 ,  0.159) \\
               &  $\widehat{\lambda}_{res}$    &    &   & (-0.031 ,  0.291) & (0.018 ,  0.202) \\
                &  $\widehat{\lambda}_{res}^{*} $    & (0.015 , 2.183)   & (0.020 , 1.271)  & (0.028 ,  0.551) & (0.044 ,  0.265) \\
\hline
\multirow{4}{*}{$(0.5, 0.5)$}  &  $\widehat{\nu}_{ls}$   & & & (0.336 ,  0.704) & (0.427 ,  0.580)   \\
                &  $\widehat{\nu}_{res}$    & (0.319 , 0.749)   & (0.366 ,  0.664)   & (0.424 ,  0.586) & (0.464 ,  0.536) \\ \cline{2-6}
                &  $\widehat{\lambda}_{ls}$    &   &  & (-0.261 ,  1.881) & (0.151 ,  0.968) \\
               &  $\widehat{\lambda}_{res}$    &    &   & (-0.103 , 1.389) & (0.121 ,  0.978) \\
                &  $\widehat{\lambda}_{res}^{*} $    & (0.099 , 9.154)   &(0.113 , 5.261)  &  (0.164 ,  2.542) & (0.241 ,  1.246) \\
  \hline
\multirow{4}{*}{$(0.75, 0.25)$}  &  $\widehat{\nu}_{ls}$   &   &   & (0.534 ,  1.019) & (0.648 ,  0.855)   \\
                &  $\widehat{\nu}_{res}$    & (0.527 , 1.057)   &(0.573 , 0.953)  & (0.648 ,  0.857) & (0.704 ,  0.798) \\ \cline{2-6}
                &  $\widehat{\lambda}_{ls}$    &    &  & (0.062 ,  0.543) & (0.117 ,  0.409) \\
               &  $\widehat{\lambda}_{res}$    &    &   & (-0.008 ,  0.618) & (0.076 ,  0.453) \\
                &  $\widehat{\lambda}_{res}^{*} $    & (0.054 , 3.070)   &(0.067 , 2.206)  & (0.087 ,  1.055) & (0.127 ,  0.573) \\
  \hline
\multirow{4}{*}{$(0.95, 5)$}  &  $\widehat{\nu}_{ls}$   &   &  & (0.700 ,  1.219) & (0.839 ,  1.067)   \\
                &  $\widehat{\nu}_{res}$    & (0.719 , 1.221)   &(0.779 , 1.150)  & (0.848 ,  1.059 & (0.904 ,  0.999) \\ \cline{2-6}
                &  $\widehat{\lambda}_{ls}$    &    &  & (-2.5731 ,  17.648) & (0.963 ,  10.334) \\
               &  $\widehat{\lambda}_{res}$    &    &   & (0.285 ,  10.951) & (1.976 ,  8.557) \\
                &  $\widehat{\lambda}_{res}^{*} $    & (1.573 , 49.872)   &(1.602 , 30.456)  & (2.006 ,  16.841) & (2.714 ,  10.067) \\
  \hline
\end{tabular*}
}
  \label{t2}
\end{table*}
\end{center}

The corresponding coverage probabilities and the widths of the interval estimates above with  a confidence level of 95\% are  displayed in Table \ref{t3}. When the sample size $n=15$, the residual-based interval estimators have minimum coverage of 90.1\% and 91.1\% for $\nu=0.95$ and $\lambda=0.1$, respectively. When $n=500$, the bootstrap interval estimator of $\lambda$ has coverage probabilities which are closer to the true confidence level than the LS-based procedure for large values of $\lambda$. However, the LS-based estimator of $\lambda$ has a better coverage than the bootstrapped residual-based  interval estimator for small $\lambda$ values. The residual-based interval estimator for $\lambda$ seemed to have slower convergence than the LS-based method. Furthermore, the residual-based estimator of $\nu$ outperformed the LS-based method as its coverage probabilities are closer to 95\%, and has narrower intervals. Overall, the coverage probabilities and interval widths still provide good merits for our estimators even when the sample size is as small as $n=15$.

\begin{center}
\begin{table*}[h!t!b!p!]
\caption{\emph{Coverage probabilities and mean widths of 95\% interval estimates  for different values of $\nu$ and  $\lambda$.}}
 \centerline {
\begin{tabular*}{6.4in}{@{\extracolsep{\fill}}c|c@{\hspace{0.1in}}||cc|cc|cc|cc}
 \hline
\multirow{2}{*}{$(\nu, \lambda)$} & \multirow{2}{*}{Estimator}  &  \multicolumn{2}{c|}{$n=15$}  &  \multicolumn{2}{c|}{$n=30$} &  \multicolumn{2}{c|}{$n=100$} & \multicolumn{2}{c}{$n=500$}  \\
 &   & Coverage & Width &  Coverage & Width & Coverage & Width & Coverage & Width  \\
  \hline \hline
\multirow{4}{*}{$(0.1, 1)$}  &  $\widehat{\nu}_{ls}$     &  &   & & &  0.949 & 0.078  &  0.952 & 0.033 \\
               &  $\widehat{\nu}_{res}$     & 0.956  & 0.098    &0.958  & 0.066  & 0.956 & 0.035  &  0.957 & 0.016  \\ \cline{2-10}
               &  $\widehat{\lambda}_{ls}$     &  &  & &  & 0.893 &  94.671 &  0.915  & 2.120  \\
               &  $\widehat{\lambda}_{res}$     &  &  & &  & 0.900 & 3.248 &  0.925 & 1.854 \\
 &  $\widehat{\lambda}_{res}^{*} $     & 0.919  & 65.411    & 0.928  & 16.424 &  0.946 & 5.131 &  0.946 & 2.161 \\
\hline
\multirow{4}{*}{$(0.25, 0.1)$}  &   $\widehat{\nu}_{ls}$     &  &    &   &  &  0.937 & 0.197 &  0.958 &  0.081  \\
               &  $\widehat{\nu}_{res}$     & 0.947  & 0.238    &0.953  & 0.162  & 0.950 & 0.087  &    0.956 & 0.038  \\ \cline{2-10}
               &  $\widehat{\lambda}_{ls}$     &   &     &  &  & 0.940 & 0.167 &  0.949 & 0.107  \\
               &  $\widehat{\lambda}_{res}$     &  &  & &  & 0.895 & 0.322 &  0.925 & 0.184   \\
               &  $\widehat{\lambda}_{res}^{*} $     & 0.911  & 2.168   & 0.938 & 1.251 &  0.952 & 0.521 &  0.951 & 0.219  \\
\hline
\multirow{4}{*}{$(0.5, 0.5)$}  &   $\widehat{\nu}_{ls}$     &   &   &  &  &  0.946 & 0.367  &  0.933 & 0.153  \\
               &  $\widehat{\nu}_{res}$     & 0.952  & 0.430    & 0.953  & 0.298  & 0.954 & 0.161&  0.949 & 0.072   \\ \cline{2-10}
               &  $\widehat{\lambda}_{ls}$     &  &    & &  & 0.931 & 2.141 &  0.923 & 0.817  \\
               &  $\widehat{\lambda}_{res}$     &  &  & &  & 0.920 & 1.492 &  0.921 & 0.857   \\
               &  $\widehat{\lambda}_{res}^{*} $     & 0.946  & 9.055    &0.950  & 5.148 &  0.947 & 2.379  &  0.948 & 1.004  \\
\hline
\multirow{4}{*}{$(0.75, 0.25)$}  &   $\widehat{\nu}_{ls}$     &   & &  & &  0.931 & 0.485 &  0.942 & 0.207  \\
               &  $\widehat{\nu}_{res}$     & 0.921  & 0.529    &0.945  & 0.379  & 0.934 & 0.209 &  0.956 & 0.094  \\ \cline{2-10}
               &  $\widehat{\lambda}_{ls}$     &   &     &  &   & 0.941 & 0.481 &  0.950 & 0.291  \\
               &  $\widehat{\lambda}_{res}$     &  &  & &  & 0.903 & 0.627 &  0.934 &  0.377   \\
               &  $\widehat{\lambda}_{res}^{*} $     & 0.916  & 3.015    &0.931  & 2.139 &  0.951 & 0.956 &  0.952 & 0.442  \\
  \hline
\multirow{4}{*}{$(0.95, 5)$}  &   $\widehat{\nu}_{ls}$     & &    &  &  &  0.923 & 0.519 &  0.945 & 0.228  \\
               &  $\widehat{\nu}_{res}$     & 0.901  & 0.502    &0.927  & 0.317  & 0.941 & 0.210 &  0.947 & 0.095  \\ \cline{2-10}
               &  $\widehat{\lambda}_{ls}$     &  &    & &  & 0.899 & 20.221 &  0.939 & 9.371  \\
               &  $\widehat{\lambda}_{res}$     &  &  & &  & 0.943 & 10.666 &  0.944 & 6.581   \\
               &  $\widehat{\lambda}_{res}^{*} $     & 0.919  & 48.299   & 0.945 & 28.854 &  0.951 & 14.748 &  0.948 & 7.356  \\
  \hline
\end{tabular*}
}
\label{t3}
\end{table*}
\end{center}

Collectively, Tables \ref{t1}--\ref{t3} strongly indicate that  the proposed point and interval  estimators performed well in our computational tests. We emphasize that the point estimates could also be regarded as reasonable starting values for better iterative estimation algorithms.

\section{Application}

We now apply our proposed methods to a real dataset. In particular, we estimate the parameters of the fractional Yule model using the branching times for plethodontid salamander dataset from  \cite{hal79} \citep[see also][]{net94a,nee01}.  The 25 data points are the times measured from each node to the present of a phylogenetic tree, and can be downloaded from the package \texttt{laser} of the R software.  The summary statistics of the inter-branching times of the plethodontid dataset are given  in Table \ref{t4} below.

\begin{center}
\begin{table*}[h!t!b!p!]
\caption{\emph{Summary statistics for the  plethodontid dataset.}}
 \centerline {
\begin{tabular*}{6.4in}{@{\extracolsep{\fill}}c|c@{\hspace{0.1in}}|c|c|c|c|c}
 \hline
Minimum & First quartile & Median & Mean & Third quartile &  Maximum & Standard deviation  \\
  \hline \hline
& & & & & & \\
0.090& 0.607 & 1.315& 3.981 & 4.405 & 22.120 & 5.966 \\
& & & & & & \\
\hline
\end{tabular*}
}
\label{t4}
\end{table*}
\end{center}

The point and the 95\%  confidence interval estimates are given in Table \ref{t5}. The LS-based point estimate (0.749) of the fractional parameter $\nu$  seemed to suggest that the plethodontid salamandar branching process is not a standard Yule process  while the residual-based point estimate (1.119) appeared to suggest otherwise.  Moreover, both the LS- and residual-based interval estimates of $\nu$  indicated that   $\nu$ could be strictly less than one, which implies  that a non-standard Yule process could model the  plethodontid salamandar dataset with a confidence level of 95\%. The  residual-based point estimate (0.011)  of $\lambda$ is more conservative  than the bootstrap- and LS-based estimate (0.049).   A similar observation can be gleaned from the 95\% interval estimates, i.e., the residual-based 95\% interval estimate is narrower than the bootstrap- and LS-based interval estimates.  

\begin{center}
\begin{table*}[h!t!b!p!]
\caption{\emph{Point and 95\% interval estimates  for  $\nu$ and  $\lambda$ of the plethodontid salamander data.}}
 \centerline {
\begin{tabular*}{3in}{c||c|c}
 \hline
Estimator &  Point estimate &  Interval estimate \\
  \hline \hline  
& & \\
$\widehat{\nu}_{ls}$ &  0.749  &  ( 0.182 ,  1.317 )    \\
& & \\
                 $\widehat{\nu}_{res}$ & 1.119    & (0.955 ,  1.283 )   \\ 
& & \\ \hline
& & \\
                 $\widehat{\lambda}_{ls}$  & 0.049   & ( 0.008 ,  0.089 )  \\
& & \\
                 $\widehat{\lambda}_{res}$     & 0.011 & ( - 0.005 , 0.027 ) \\
& & \\
   $\widehat{\lambda}_{res}^{*} $     &   &  ( 0.003 , 0.051 ) \\
& & \\
\hline
\end{tabular*}
}
\label{t5}
\end{table*}
\end{center}

We also tested the residuals for normality using the Shapiro-Wilk,    Anderson-Darling, Cramer-von Mises, Lilliefors, Pearson chi-square, and the  Shapiro-Francia tests, which gave the $p$-values  0.811, 0.651, 0.619, 0.609, 0.849, and  0.461, correspondingly. Hence,  these  $p$-values indicated good fit of the  fractional Yule process to the plethodontid salamandar  data.

\section{Concluding remarks}

We have proposed  closed-form expressions of the estimators of the parameters $\nu$ and $\lambda$ for the fractional  linear birth or Yule, the  fractional linear  death,  and the fractional sublinear death processes.  The estimators were derived by taking advantage of the known structural form of the logarithm of the  random inter-event times and the well-studied  least squares regression procedure. The explicit formulas led to computationally simple and fast parameter estimation  procedures. The inter-death time distributions and variances  of  the fractional linear and sublinear death processes  were also obtained. These statistical properties were necessary for generating sample trajectories and for our estimation procedures to be applicable in these  processes. It has also been shown  that the proposed procedure can be easily extended to certain models that have different model parameterizations than the linear ones. The proposed methods were used to model a real physical process. Generally, the extensive computational tests showed favorable results for the proposed estimators.

We cite some extensions which would be worth pursuing in the future. For instance, improving the small-sample performance of the least squares-based estimators and  developing other estimators using the likelihood  approach or a  re-sampling technique  would be valuable pursuits. The application of these methods in practice,  and the characterization of the appropriate functions $m(\theta)$ and $q(j)$ would also be of interest.


\bibliographystyle{spbasic}      

\begin{thebibliography}{}


\bibitem[Aldous (2001)]{dja01} Aldous, D.J.:   Stochastic models and descriptive statistics for phylogenetic trees,
from Yule to today, Statistical Science 16(1), 23--34 (2001)


\bibitem[Cahoy and Polito (2012)]{cap12} Cahoy D.O., Polito, F.:  Simulation and estimation for the fractional Yule process. Methodology and Computing in App. Prob. 14(2),  383--403  (2012) 

\bibitem[Cahoy et al. (2010)]{cuw10} Cahoy, D.O., Uchaikin, V.V., Woyczynski, W.A.: Parameter estimation for fractional Poisson processes. J. Stat. Plan. Inf. 140, 3106--3120  (2010)


\bibitem[Efron and Tibshirani (1994)]{eat94} Efron, B., Tibshirani, R.:  An Introduction to the Bootstrap. CRC Press, Boca Raton, FL (1994)


\bibitem[Ferguson (1996)]{fer96} Ferguson, T.: A Course in Large Sample Theory, Chapman \& Hall, Great Britain (1996)

\bibitem[Highton and Larson (1979)]{hal79} Highton, R.,  Larson,  L. A.:  The genetic relationships of
the salamanders of the genus Plethodon. Syst. Zool. 28, 579--599 (1979)


\bibitem[Montgomery et al. (2006)]{mpv06} Montgomery, D.C.,  Peck, E.A., Vining, G.G.: Introduction to linear regression analysis 4 ed, John Wiley \&  Sons, Inc., Great Britain, United States of America (2006)


\bibitem[Nee et al. (1994a)]{net94a}Nee, S., Holmes, E.C.,  May, R.M., Harvey, P.H.:  Extinction rates can be estimated from molecular phylogenies. Philos. Trans. R. So. Lond. B 344, 77-82 (1994a)

\bibitem[Nee (2001)]{nee01}Nee, S.:  Inferring speciation rates from phylogenies. Evolution 55, 661-668. (2001)

\bibitem[Orsingher and Polito (2010)]{oap10} Orsingher, E., Polito, F.: Fractional pure birth processes. Bernoulli 16, 858--881  (2010)

\bibitem[Orsingher and Polito (2011)]{oap11} Orsingher, E., Polito, F.: Randomly stopped nonlinear fractional birth processes. Submitted  (2011)

\bibitem[Orsingher et al. (2010)]{ops10} Orsingher, E., Polito, F., Sakhno, L.: Fractional non-Linear, linear and sublinear death processes. J.  Stat. Phys. 141, 68--93. (2010).

\bibitem[Paradis (2012)]{pe12}Paradis, E.:  Analysis of Phylogenetics and Evolution with R, 2 ed, Springer, New York, USA (2012)


\bibitem[Uchaikin et al. (2008)]{ucs08} Uchaikin, V.V., Cahoy, D.O., Sibatov, R.T.: Fractional processes: from Poisson to branching one. International J. Bifurcation 18, 2717--2725  (2008)

\end{thebibliography}

\end{document}